\documentclass[12pt,reqno]{amsart}
\usepackage{amssymb,verbatim,enumerate,ifthen,graphics,graphicx}
\usepackage[mathscr]{eucal}
\usepackage[utf8]{inputenc}
\usepackage[T1]{fontenc}
\newcommand{\R}{\mathbb{R}}
\newcommand{\N}{\mathbb{N}}
\newcommand{\Z}{\mathbb{Z}}
\newcommand{\D}{\mathbb{D}}
\renewcommand{\P}{\mathscr{P}}
\newcommand{\U}{\mathscr{U}}

\newcommand{\dist}{\mathop{\hbox{\rm dist}}\nolimits}

\newcommand{\rec}{\mathop{\hbox{\rm rec}}\nolimits}

\newcommand{\cl}{\mathop{\hbox{\rm cl}}\nolimits}

\newtheorem{theorem}{Theorem}

\newtheorem*{theorem*}{Theorem}
\def\Thm#1#2{\ifthenelse{\equal{#1}{*}}{\begin{theorem*}#2\end{theorem*}}
  {\begin{theorem}\label{T#1}#2\end{theorem}}}
\newtheorem{Atheorem}{Theorem}

\def\thm#1{Theorem~\ref{T#1}}
\newtheorem{proposition}[theorem]{Proposition}
\newtheorem*{proposition*}{Proposition}
\def\Prp#1#2{\ifthenelse{\equal{#1}{*}}{\begin{proposition*}#2\end{proposition*}}
             {\begin{proposition}\label{P#1}#2\end{proposition}}}
\def\prp#1{Proposition~\ref{P#1}}
\newtheorem{corollary}[theorem]{Corollary}
\newtheorem*{corollary*}{Corollary}
\def\Cor#1#2{\ifthenelse{\equal{#1}{*}}{\begin{corollary*}#2\end{corollary*}}
             {\begin{corollary}\label{C#1}#2\end{corollary}}}
\def\cor#1{Corollary~\ref{C#1}}
\newtheorem{lemma}[theorem]{Lemma}
\newtheorem*{lemma*}{Lemma}
\def\Lem#1#2{\ifthenelse{\equal{#1}{*}}{\begin{lemma*}#2\end{lemma*}}
             {\begin{lemma}\label{L#1}#2\end{lemma}}}
\def\lem#1{Lemma~\ref{L#1}}
\newtheorem{example}[theorem]{Example}
\newtheorem*{example*}{Example}
\long\def\Exa#1#2{\ifthenelse{\equal{#1}{*}}{\begin{example*}\rm #2\end{example*}}
            {\begin{example}\label{Ex#1}\rm #2\end{example}}}

\newtheorem{problem}[subsection]{Problem}

\theoremstyle{definition}
\newtheorem{definition}[theorem]{Definition}
\newtheorem*{definition*}{Definition}
\def\Defi#1#2{\ifthenelse{\equal{#1}{*}}{\begin{definition*}#2\end{definition*}}
      {\begin{definition}\label{D#1}#2\end{definition}}}

\newtheorem{remark}[theorem]{Remark}
\newtheorem*{remark*}{Remark}
\long\def\Rem#1#2{\ifthenelse{\equal{#1}{*}}{\begin{remark*}#2\end{remark*}}
             {\begin{remark}\label{R#1}#2\end{remark}}}

\def\eq#1{{\rm(\ref{E#1})}}
\def\Eq#1#2{\ifthenelse{\equal{#1}{*}}
  {\begin{equation*}\begin{aligned}[]#2\end{aligned}\end{equation*}}
  {\begin{equation}\begin{aligned}[]\label{E#1}#2\end{aligned}\end{equation}}}

\allowdisplaybreaks

\begin{document}

\begin{flushright}
\end{flushright}
\vspace{5mm}

\date{\today}

\title[Bernstein--Doetsch type theorems with Tabor type error terms]
{Bernstein--Doetsch type theorems \\ with Tabor type error terms for set-valued maps}

\author[A. Gil\'anyi]{Attila Gil\'anyi}
\address{Faculty of Informatics, University of Debrecen, H-4032 Debrecen, Egyetem t\'er 1, Hungary}
\email{gilanyi@inf.unideb.hu}

\author[C. Gonz\'alez]{Carlos Gonz\'alez}
\address{Escuela de Matematicas, Universidad Central de Venezuela, Caracas, Venezuela}
\email{carlosl.gonzalez@ciens.ucv.ve}

\author[K. Nikodem]{Kazimierz Nikodem}
\address{Department of Mathematics, University of Bielsko-Biala, 
ul.\ Willowa~2, 43-309 Bielsko-Bia{\l}a, Poland}
\email{knikodem@ath.bielsko.pl}

\author[Zs. P\'ales]{Zsolt P\'ales}
\address{Institute of Mathematics, University of Debrecen, H-4032 Debrecen, Egyetem t\'er 1, Hungary}
\email{pales@science.unideb.hu}

\dedicatory{Dedicated to the 100th anniversary of the Bernstein--Doetsch Theorem}


\thanks{This research was supported by the Hungarian Scientific Research Fund (OTKA) Grant K 111651.}

\begin{abstract}
In this paper, we investigate set-valued maps of strongly and approximately Jensen convex and Jensen 
concave type. We present counterparts of the Bernstein--Doetsch Theorem with Tabor type error terms.
\end{abstract}

\maketitle

\section{Introduction}

In the theory of convex functions the celebrated Bernstein--Doetsch Theorem published exactly hundred years ago 
(\cite{BerDoe15}) plays a fundamental role. Roughly speaking this theorem states that if a function is Jensen convex 
and possesses the local upper boundedness property then it is also convex. This result has been extended in several 
ways in several directions. These extensions concern vector-valued and set-valued mappings and also the approximately 
and strongly convex setting. Such extensions were obtained, among others, in the papers \cite{AveCar90},
\cite{AzoGimNikSan11}, \cite{CarNikPap93}--\cite{Tru84}. 
\nocite{GonNikPalRoa14,Haz05a,Haz07b,HazPal04,HazPal05,HazPal09,Lac99,LeiMerNikSan13,MakPal10b,MakPal12b,%
MakPal12c,MakPal13b,MurTabTab12,NgNik93,Nik86,Nik87a,Nik87c,Nik89,Pap90,TabTab09a,TabTab09b,TabTabZol10a,%
TabTabZol10b,Tru84}

A summary and detailed comparison of these generalizations can be found in the recent paper \cite{GonNikPalRoa14}.

To recall the two main results established in this paper, let $X$ and $Y$ be Hausdorff topological linear spaces. 
Denote by $\P_0(Y)$ the family of nonempty subsets of $Y$ and let $D\subseteq X$ be a nonempty convex set and assume 
that $A,B:(D-D)\to\P_0(Y)$ are set-valued maps such that $0\in A(x)\cap B(x)$ for all $x\in (D-D)$. 

In terms of these data, for a set-valued 
map $F:D\to\P_0(Y)$, the following Jensen convexity and Jensen concavity type inclusions can be considered:
\Eq{JCV0}{
\dfrac{F(x) + F(y)}{2} + A(x-y) \subseteq \cl\bigg(F\bigg(\dfrac{x+y}{2}\bigg) + B(x-y)\bigg) 
   \qquad (x,y\in D)
}
and
\Eq{JCC0}{
F\bigg(\dfrac{x+y}{2}\bigg) + A(x-y) \subseteq \cl\bigg(\dfrac{F(x) + F(y)}{2} + B(x-y)\bigg) \qquad (x,y\in D).
}
If, in particular, $A=\{0\}$ or $B=\{0\}$ then we call the above properties approximate and strong Jensen 
convexity/concavity, respectively.

The basic problems addressed in \cite{GonNikPalRoa14} were to derive convexity and concavity type inclusions for $F$ 
under certain additional regularity properties. 

Given a set-valued map $S:(D-D)\to\P_0(Y)$, define the Takagi transformation $S^T:\R\times D\to\P_0(Y)$ by
\Eq{Tak0}{
  S^T(t,u):=\cl\bigg(\bigcup_{n=0}^{\infty} \sum_{k=0}^{n} 
                 \frac{1}{2^k}S\big(2d_{\Z}(2^kt)u\big)\bigg)\qquad(t\in\R,\,u\in D-D),
}
where $d_{\Z}:\R\to\R$ is defined by
\Eq{dz}{
  d_{\Z}(t):=\dist(t,\Z):=\inf_{k\in\Z}|k-t| \qquad(t\in\R).
}
Denote by $K$ the closure of the recession cone of the set-valued map $B$ (for its definition, see the next section).

The main results of the paper \cite{GonNikPalRoa14} can be formulated as follows.

\Thm{Convex}{Assume that the set-valued map $F$ satisfies the Jensen convexity type inclusion \eq{JCV0} 
and it possesses the following properties:
\begin{enumerate}[(i)]
 \item $F$ is pointwise closedly $K$-lower bounded, i.e., for each $x\in D$, there exists a bounded set 
$H\subseteq Y$ such that $F(x)\subseteq \cl(H+K)$;
 \item $F$ is locally closedly weakly $K$-upper bounded on $D$, i.e., for all $x\in D$, there exists
an open set $U$ containing $x$ and a bounded set $H\subseteq Y$ such that 
$0\in \cl(F(u)+H+K)$ holds for all $u\in D\cap U$.
\end{enumerate}
Then $F$ satisfies the convexity type inclusion
\Eq{*}{
 tF(x)+(1-t)F(y)+A^T(t,x-y)\subseteq \cl\big(F(tx+(1-t)y)+B^T(t,x-y)\big) \\
   \quad(x,y\in D,\,t\in[0,1]). 
}}

\Thm{Concave}{Assume that the set-valued map $F$ satisfies the Jensen concavity type inclusion \eq{JCC0} 
and it possesses the following properties:
\begin{enumerate}[(i)]
 \item $F$ is pointwise closedly $K$-convex, i.e., $tF(x)+(1-t)F(x)\subseteq \cl(F(x) + K)$ holds 
for each $x\in D$ and for all $t\in[0,1]$;
 \item $F$ is locally closedly $K$-lower bounded, i.e., for each $x\in D$, there exists a neighborhood
$U$ of $x$ and a bounded set $H\subseteq Y$ such that $F(u)\subseteq \cl(H+K)$ for all $u\in D\cap U$.
\end{enumerate}
Then $F$ satisfies the concavity type inclusion
\Eq{*}{
 F(tx+(1-t)y)+A^T(t,x-y)\subseteq \cl\big(tF(x)+(1-t)F(y)+B^T(t,x-y)\big) \\
  \quad(x,y\in D,\,t\in[0,1]). 
}}

The above so-called Takagi type error term for the first time appeared in the paper \cite{HazPal04} by H\'azy and 
P\'ales. Its sharpness was proved by Boros \cite{Bor08} and, in a more general setting, by Mak\'o and P\'ales
\cite{MakPal10b}. A new type of error term has been discovered by Jacek Tabor and J\'ozef Tabor in \cite{TabTab09a} and 
its sharpness has been investigated by Mak\'o and P\'ales \cite{MakPal13b} in the real valued setting.  Motivated by 
these results, in this paper we introduce the set-valued and vector-valued generalization of the Tabor type error term 
and deduce results that are completely analogous to \thm{Convex} and \thm{Concave} but involve the Tabor type error 
term instead of the Takagi type one.

\section{Terminology and Auxiliary Results}
\setcounter{theorem}{0}

Let $X$ be a linear space and $Y$ be a Hausdorff topological vector space in the sequel.
A base of balanced neighborhoods of zero in $Y$ will be denoted by $\U(Y)$. The notation $\P_0(Y)$ will stand 
for the class of nonempty subsets of $Y$. For the sake of brevity, for a set $S\subseteq Y$ and
$n\in\N$, we set
\Eq{*}{
  [n]S:=\{x_1+\cdots+x_n\mid x_1,\dots,x_n\in S\}.
}

Given a convex cone $K\subseteq Y$, we recall from \cite{GonNikPalRoa14} and introduce several 
notions that will be needed for the formulation of our results.

A set $S\subseteq Y$ is called \emph{closedly $K$-convex} if, for all $t\in[0,1]$,
\Eq{*}{
  tS+(1-t)S\subseteq \cl(S+K).
}
For a closedly $K$-convex set $S$, it easily follows by an inductive argument that, for all $t_1,\dots,t_n\geq0$ with 
$t_1+\cdots+t_n=1$, we have
\Eq{*}{
  t_1S+\cdots+t_nS\subseteq \cl(S+K)
}
and hence the convex hull of $S$ is contained in $\cl(S+K)$.

We say that $S\subseteq Y$ is \emph{closedly $K$-starshaped with respect to $y\in Y$} if,
for all $t\in[0,1]$,
\Eq{*}{
  tS+(1-t)y\subseteq \cl(S+K).
}
By putting $t=0$, the above inclusion implies that $y\in\cl(S+K)$ provided that $S$ is closedly $K$-starshaped with 
respect to $y$. If $S$ is closedly $K$-starshaped with respect to $0$, then we will briefly say that $S$ is 
\textit{closedly $K$-starshaped}. Obviously, convex sets are closedly $K$-convex, starshaped sets are closedly 
$K$-starshaped and closedly $K$-convex sets are closedly $K$-starshaped with respect to their elements.
In the case when $K=\{0\}$ we simply omit $K$ and call $S$ closedly convex and closedly starshaped, respectively.

A set $H\subseteq Y$ is called \emph{bounded} if, for all $U\in\U(Y)$, there exists a positive number
$t$ such that $tH\subseteq U$. One can easily see that the family of bounded sets is closed under
finite union, under algebraic addition and multiplication by scalars. 

A set $S\subseteq Y$ is called \emph{closedly $K$-lower bounded}, if there exists a bounded set 
$H\subseteq Y$ such that $S\subseteq\cl(H+K)$. One can prove that this class of sets
is also closed under finite union, under algebraic addition and multiplication by scalars.

For a nonempty subset $S\subseteq Y$, we define its \textit{recession cone} by
\Eq{*}{
  \rec(S):=\{y\in Y\mid ty+S\subseteq S\mbox{ for all }t\geq0\}.
}
The \textit{recession cone of a set-valued map} $S:D\to\P_0(Y)$ is introduced by
\Eq{*}{
 \rec(S):=\bigcap_{x\in D} \rec S(x).
}

The basic properties of the recession cone are summarized in the following lemma
which was established in the paper \cite{GonNikPalRoa14}.

\Lem{rec}{Let $S\subseteq Y$ be a nonempty set. Then
\begin{enumerate}[(i)]
 \item $\rec(S)$ is a convex cone containing $0$;
 \item $K=\rec(S)$ is the largest cone $K$ such that $K+S\subseteq S$ is valid;
 \item $\overline{\rec}(S)\subseteq\rec(\overline{S})$;
 \item for all $y\in Y$, $t>0$, $\rec(y+tS)=\rec(S)$;
 \item for all nonempty sets $S,T\subseteq Y$, $\rec(S)+\rec(T)\subseteq \rec(S+T)$.
\end{enumerate}}

Let $F:D\to\P_0(Y)$ be a set-valued map. We say that $F$ is
\emph{directionally $K$-upper semicontinuous at a point $p\in D$} if, for every direction $h\in X$ 
and for every neighborhood $U\in\U(Y)$, there exists a positive number $\delta$ such that
\Eq{*}{
  F(p+th)\subseteq F(p)+U+K
}
for all $0<t<\delta$ with $p+th\in D$. In the particular case when $K=\{0\}$, we speak about \textit{directional
upper semicountinuity at $p$}. Clearly, directionally upper semicontinuous maps are also directionally $K$-upper 
semicontinuous for any cone $K$.

Analogously, we say that $F$ is
\emph{directionally $K$-lower semicontinuous at a point $p\in D$} if, for every direction $h\in X$ 
and for every neighborhood $U\in\U(Y)$, there exists a positive number $\delta$ such that
\Eq{*}{
  F(p)\subseteq F(p+th)+U+K
}
for all $0<t<\delta$ with $p+th\in D$. If $F$ is simultaneously directionally $K$-upper and $K$-lower semicontinuous
at $p$, then we call $F$ \emph{directionally $K$-continuous at $p$}.

\Lem{dusc}{Let $K\subseteq Y$ be a convex cone and $S,T\subseteq Y$ be nonempty and closedly $K$-lower bounded sets 
which are closedly $K$-starshaped with respect to some elements of $Y$. Then, the set-valued 
mapping $t\mapsto tS+(1-t)T$ is directionally $K$-continuous on $[0,1]$.}

\begin{proof}
Let $U\in\U(X)$. It is sufficient to show that there exists a positive real number $\delta$ such that, for 
$s,t\in[0,1]$ with $|t-s|<\delta$, 
\Eq{st}{
  tS+(1-t)T\subseteq sS+(1-s)T + U+K.
}

First choose $V\in\U(Y)$ such that $[3]V\subseteq U$. Obviously, \eq{st} holds for $s=t$. Without loss of generality, 
we may assume that $0\leq s<t\leq 1$. Suppose that $T$ is closedly $K$-starshaped with respect to $v\in Y$. Then we have
\Eq{est}{
tS+(1-t)T
&\subseteq sS + (t-s)S + (1-t)T + (t-s)v+(s-t)v \\
&= sS + (t-s)S + (1-s)\bigg(\frac{1-t}{1-s}T+ \frac{t-s}{1-s}v\bigg)+(s-t)v \\
&\subseteq sS + (t-s)(S-v) + (1-s)\cl(T+K) \\
&\subseteq sS + (1-s)T + (t-s)(S-v) +V+K.
}
As the set $(S-v)$ is closedly $K$-lower bounded, there exists a bounded set 
$H\subseteq Y$ such that $S-v\subseteq\cl(H+K)$. By the boundedness of the set $H$, there exists
a positive number $\delta\leq 1$ such that $\delta H\subseteq V$. 

Therefore, if $t-s<\delta$, then
\Eq{*}{
  (t-s)(S-v)\subseteq (t-s)\cl(H+K)\subseteq (t-s)(H+V+K) \subseteq (t-s)H+V+K \subseteq [2]V+K.
}
Combining this estimate with \eq{est}, we obtain that \eq{st} holds. The proof for the case
$t<s$ is completely analogous and uses that $S$ is $K$-starshaped and that $T$ is closedly $K$-lower bounded.
\end{proof}

\section{Tabor transformation of set-valued maps}
\setcounter{theorem}{0}

Assume now that $D\subseteq X$ is a starshaped set and consider a set-valued map $S:D\to\P_0(Y)$.
For such a map, we define $S^\perp:\R\times D\to Y$ by
\Eq{Tab0}{
  S^\perp(t,x):=\cl\bigg(\bigcup_{n=0}^{\infty} \sum_{k=0}^{n}
                 2d_{\Z}(2^kt)S\Big(\frac{x}{2^k}\Big)\bigg)\qquad((t,x)\in\R\times D),
}
where $d_{\Z}:\R\to\R$ is defined by \eq{dz}.
The set-valued map $S^\perp$ will be called the \textit{Tabor transformation} of the map $S$ in the sequel.
In view of the property that $d_\Z$ vanishes on $\Z$, it is immediate to see that $S^\perp(t,x)=\{0\}$ for all 
$(t,x)\in\Z\times D$.

In the following lemma we establish the relationship between a set-valued map and its Tabor transformation
and also between the recession cones of these set-valued maps.

\Lem{TT}{Let $D\subseteq X$ be a starshaped set and $S:D\to\P_0(Y)$. Then
\Eq{TT1}{
   S^\perp\big(\tfrac12,x)\big)=\cl(S(x)) \qquad (x\in D)
}
and 
\Eq{TT2}{
  \overline{\rec}\, S\subseteq \rec \big(S^\perp(t,x)\big) \qquad((t,x)\in(\R\setminus\Z)\times D).
}}

\begin{proof} Observe that $d_{\Z}\big(\tfrac12\big)=\tfrac12$ and $d_{\Z}\big(2^k\cdot\tfrac12\big)=0$ for $k\in\N$.
Thus,
\Eq{*}{
  S^\perp\big(\tfrac12,x\big)
   =\cl\bigg(\bigcup_{n=0}^{\infty} \sum_{k=0}^{n}2d_{\Z}\Big(\frac{2^k}2\Big)S\Big(\frac{x}{2^k}\Big)\bigg)
   =\cl(S(x))
}
for all $x\in D$.

For the proof of \eq{TT2}, denote $K:=\rec S$. Then, for all $x\in D$, we have 
\Eq{*}{
   S(x)+K\subseteq S(x).
}
Let $(t,x)\in(\R\setminus\Z)\times D$ be arbitrary. Using the above inclusion and $d_Z(t)\neq0$, for all $n\geq0$,
we get
\Eq{*}{
  \bigg(\sum_{k=0}^{n}2d_{\Z}(2^kt)S\Big(\frac{x}{2^k}\Big)\bigg)+K
  =\sum_{k=0}^{n}2d_{\Z}(2^kt)\Big(S\Big(\frac{x}{2^k}\Big)+K\Big)
  \subseteq \sum_{k=0}^{n}2d_{\Z}(2^kt)S\Big(\frac{x}{2^k}\Big)\subseteq S^\perp(t,x).
}
Therefore,
\Eq{*}{
  \bigcup_{n=0}^\infty \bigg(\sum_{k=0}^{n}2d_{\Z}(2^kt)S\Big(\frac{x}{2^k}\Big)\bigg)+K
  \subseteq S^\perp(t,x),
}
which implies
\Eq{*}{
  S^\perp(t,x)+\cl(K)\subseteq S^\perp(t,x).
}
This completes the proof of \eq{TT2}.
\end{proof}

The following result implies that the Tabor transformation of a set-valued map which is constructed as the product of a 
nonnegative scalar function and a convex subset of $Y$ is the product of the Tabor transformation of the scalar 
function and of the same set.

\Prp{Tab}{Let $D\subseteq X$ be a starshaped set, $K\subseteq Y$ be a convex cone and $S_0\subseteq Y$ be a closedly 
$K$-convex and closedly $K$-starshaped set. Let $\varphi:D\to \R_+$ be a nonnegative function such that
\Eq{phi}{
  \sum_{k=0}^\infty \varphi\Big(\frac{x}{2^n}\Big)<\infty \qquad(x\in D).
}
Define $S:D\to\P_0(Y)$ by $S(x):=\varphi(x)S_0+K$. Then
\Eq{Tab1}{
  S^\perp(t,x)=\cl\big(\varphi^\perp(t,x) S_0+K\big) \qquad ((t,x)\in(\R\setminus\Z)\times D), 
}
where
\Eq{Tab2}{
  \varphi^\perp(t,x):=\sum_{n=0}^{\infty}2d_{\Z}(2^nt)\varphi\Big(\frac{x}{2^n}\Big)
  \qquad((t,x)\in\R\times D).
}}

\begin{proof} Using the convergence condition \eq{phi} and the estimate $0\leq 2d_\Z\leq1$, it follows that the series 
in \eq{Tab2} converges uniformly in $t$ and consequently, for all $x\in D$, the mapping $t\mapsto\varphi^\perp(t,x)$ is 
continuous on $\R$.

To prove \eq{Tab1}, fix $(t,x)\in(\R\setminus\Z)\times D$. Then $d_\Z(t)>0$.
If, for an element $x\in D$, we have $\varphi^\perp(t,x)=0$, then (by the nonnegativity of $\varphi$) it follows that
$d_{\Z}(2^nt)\varphi\big(\frac{x}{2^n}\big)=0$ for all $n\geq0$. Then,
\Eq{*}{
  S^\perp(t,x)
  &=\cl\bigg(\bigcup_{n=0}^{\infty} \sum_{k=0}^{n}2d_{\Z}(2^kt)\Big(\varphi\Big(\frac{x}{2^k}\Big)S_0+K\Big)\bigg)\\
  &=\cl\bigg(\bigcup_{n=0}^{\infty} 
      \Big(\Big(\sum_{k=0}^{n}2d_{\Z}(2^kt)\varphi\Big(\frac{x}{2^k}\Big)S_0\Big)+K\Big)\bigg)
  =\cl(K)=\cl\big(\varphi^\perp(t,x) S_0+K\big),
}
which shows the validity of \eq{Tab1} in this case. 
Therefore, in the sequel, we may also assume that $x\in D$ is chosen such that $\varphi^\perp(t,x)\neq0$.

Using that $S_0$ is closedly $K$-convex and closedly $K$-starshaped, we have, for all $t_0,\dots,t_n\geq0$ and 
$t_0+\cdots+t_n\leq t$, that
\Eq{ttt}{
  t_0S_0+\cdots+t_nS_0\subseteq \cl(tS_0+K).
}
Indeed, if $t_1=\cdots=t_n=t=0$, then \eq{ttt} is equivalent to the trivial inclusion $0\in\cl(K)$. 
If $t_1=\cdots=t_n=0<t$, then \eq{ttt} reduces to $0\in t\cl(S_0+K)$, which is a consequence of the closed 
$K$-starshapedness of $S_0$. Finally, if $t_0+\cdots+t_n>0$, then using that $S_0$ is closedly $K$-convex and closedly 
$K$-starshaped, we obtain
\Eq{*}{
  t_0S_0+\cdots+t_nS_0
  &=(t_0+\cdots+t_n)\Big(\frac{t_0}{t_0+\cdots+t_n}S_0+\cdots+\frac{t_n}{t_0+\cdots+t_n}S_0\Big)\\
  &\subseteq(t_0+\cdots+t_n)\cl(S_0+K)=\cl((t_0+\cdots+t_n)S_0+K)\\
  &=\cl\Big(t\frac{t_0+\cdots+t_n}{t}S_0+K\Big)\subseteq \cl(t\cl(S_0+K)+K)=\cl(tS_0+K).
}

To complete the proof of \eq{Tab1} in the general case, we verify first the inclusion $\subseteq$.
Applying \eq{ttt}, for all $n\geq 0$, we obtain 
\Eq{*}{
\sum_{k=0}^{n}2d_{\Z}(2^kt)\varphi\Big(\frac{x}{2^k}\Big)S_0
\subseteq\cl\Bigg(\bigg(\sum_{k=0}^{\infty}2d_{\Z}(2^kt)\varphi\Big(\frac{x}{2^k}\Big)\bigg)S_0+K\Bigg)
  =\cl\big(\varphi^\perp(t,x)S_0+K\big).
}
Therefore (using also that $d_Z(t)>0$),
\Eq{*}{
  \sum_{k=0}^{n}2d_{\Z}(2^kt)\Big(\varphi\Big(\frac{x}{2^k}\Big)S_0+K\Big)
  &\subseteq \bigg(\sum_{k=0}^{n}2d_{\Z}(2^kt)\varphi\Big(\frac{x}{2^k}\Big)S_0\bigg)+K \\
  &\subseteq \cl\big(\varphi^\perp(t,x)S_0+K\big)+K=\cl\big(\varphi^\perp(t,x)S_0+K\big).
}
Since this last inclusion is valid for all $n\geq0$, the desired inclusion follows immediately. 

For the opposite inclusion, in view of $\varphi^\perp(t,x)\neq0$, we can choose $n_0$ such that 
$d_{\Z}(2^{n_0}t)\varphi\big(\frac{x}{2^{n_0}}\big)>0$. Now observe that, for every $n\geq n_0$,
\Eq{*}{
\bigg(\sum_{k=0}^n & 2d_\Z(2^kt)\varphi\Big(\frac{x}{2^k}\Big)\bigg)S_0 + K
  =\bigg(\sum_{k=0}^n2d_\Z(2^kt)\varphi\Big(\frac{x}{2^k}\Big)\bigg)(S_0+K) \\
  &\subseteq\sum_{k=0}^n2d_\Z(2^kt)\bigg(\varphi\Big(\frac{x}{2^k}\Big)S_0+\varphi\Big(\frac{x}{2^k}\Big)K\bigg)
  \subseteq\sum_{k=0}^n2d_\Z(2^kt)\bigg(\varphi\Big(\frac{x}{2^k}\Big)S_0+\cl(K)\bigg)\\
  &\subseteq\cl\bigg(\sum_{k=0}^n2d_\Z(2^kt)\bigg(\varphi\Big(\frac{x}{2^k}\Big)S_0+K\bigg)\bigg)
  \subseteq\cl\bigg(\bigcup_{\ell=0}^\infty\sum_{k=0}^\ell 2d_\Z(2^kt)S\Big(\frac{x}{2^k}\Big)\bigg)
  = S^\perp(t,x).
}
Let $u\in K$ and $v\in S_0$. Then, the above inclusion yields that, for every $n\geq n_0$,
\Eq{*}{
  y_n:=\sum_{k=0}^n2d_\Z(2^kt)\varphi\Big(\frac{x}{2^k}\Big)v+u\in S^\perp(t,x).
}
The right hand side is a closed set, hence the limit of the sequence $(y_n)$ is also contained in $S^\perp(t,x)$.
Therefore, for all $u\in K$ and $v\in S_0$,
\Eq{*}{
  \varphi^\perp(t,x)v+u
  =\sum_{k=0}^\infty 2d_\Z(2^kt)\varphi\Big(\frac{x}{2^k}\Big)v+u
  =\lim_{n\to\infty} y_n
  \in S^\perp(t,x).
}
which implies that
\Eq{*}{
\varphi^\perp(t,x)S_0+K\subseteq S^\perp(t,x).
}
This completes the proof of the reversed inclusion in \eq{Tab1}.
\end{proof}

\Cor{Tak}{Let $X$ be a normed space, $D\subseteq X$ be a starshaped set, $K\subseteq Y$ be a convex cone, 
$S_0\subseteq Y$ be a closedly $K$-convex set containing $0\in Y$, and $\alpha>0$. Define $S:D\to\P_0(Y)$ by
$S(x):=\|x\|^\alpha S_0+K$. Then
\Eq{Tab1+}{
  S^\perp(t,x)
  =\cl\big(\tau_\alpha(t)\|x\|^\alpha S_0+K\big) \qquad((t,x)\in(\R\setminus\Z)\times D), \\
}               
where $\tau_\alpha:\R\to\R$ is defined by 
\Eq{Tab2+}{
  \tau_\alpha(t):=\sum_{n=0}^{\infty}2^{1-\alpha n} d_{\Z}(2^nt)
  \qquad(t\in\R).
}}

\begin{proof} To obtain the statement, we can apply \prp{Tab} with the function $\varphi$
defined by $\varphi(x):=\|x\|^\alpha$. Observe that
\Eq{*}{
  \sum_{n=0}^{\infty}\varphi\Big(\frac{x}{2^n}\Big)
	=\sum_{n=0}^{\infty}\frac{\|x\|^\alpha}{2^{\alpha n}}
	=\frac{2^\alpha}{2^\alpha-1}\|x\|^\alpha < \infty
	\qquad(x\in D)
}
and
\Eq{*}{
  \varphi^\perp(t,x)=\sum_{n=0}^{\infty}2d_{\Z}(2^nt)\varphi\Big(\frac{x}{2^n}\Big)
	=\sum_{n=0}^{\infty}2d_{\Z}(2^nt)\frac{\|x\|^\alpha}{2^{\alpha n}}
	=\tau_\alpha(t)\|x\|^\alpha
	\qquad((t,x)\in\R\times D).
}
Therefore, \eq{Tab1+} is a consequence of formula \eq{Tab1} of \prp{Tab}.
\end{proof}

\Rem{Tak}{An important particular case is when $\alpha=1$, then $\tau_1(t)=2T$, 
where $T$ is the Takagi function defined by 
\Eq{takagi}{
T(t):=\sum_{n=0}^\infty\frac{d_\Z(2^nt)}{2^n} \qquad(t\in\R)
} 
Observe that $\tau_\alpha$ (for any $\alpha>0$) satisfies the functional equation
\Eq{TT}{
  \tau_\alpha(t)=2d_{\Z}(t)+\frac1{2^\alpha}\tau_\alpha(2t) \qquad(t\in\R).
}
By applying the Banach Fixed Point Theorem, this functional equation has a unique solution in the Banach
space of bounded real functions over the real line (which is equipped with the supremum norm).}

The following lemma will be useful in order to prove our main theorems related to the convexity and concavity type 
inclusion. For its formulation, we introduce a convergence property for a sequence of sets in $Y$. Given a convex cone 
$K\subseteq Y$ and a sequence $S_k\subseteq Y$, this sequence is called serially $K$-Cauchy if, for all $U\in\U(Y)$, 
there exists $m\in\N$ such that, for all $n\geq m$,
\Eq{iii}{
 \sum_{k=m}^nS_k \subseteq U+K.
}

\Lem{CchyConv}{
Let $K\subseteq Y$ be a convex cone and let $(S_k)$ be a sequence of nonempty subsets of $Y$ such that 
\begin{enumerate}[(i)]
  \item For all $k\geq 0$, the set $S_k$ is closedly $K$-starshaped and closedly $K$-lower bounded;
  \item The sequence $(S_k)$ is serially $K$-Cauchy. 
\end{enumerate}
Then, for all $U\in\U(Y)$, there exists a positive number $\delta$ such that,
for all $t,s\in\R$ with $|t-s|<\delta$, 
\Eq{ts0}{
\cl\bigg(\bigcup_{n=0}^\infty\sum_{k=0}^n d_\Z(2^kt)S_k\bigg) 
  \subseteq \cl\bigg(\bigcup_{n=0}^\infty\sum_{k=0}^n d_\Z(2^ks)S_k\bigg) + U + K.
}}

\begin{proof}
Let $U\in\U(Y)$. Then there exists $V\in\U(Y)$ such that $[6]V\subseteq U$. 

First, we are going to show that there exists a positive number $\delta$ such that,
for all $n\geq0$ and for all $t,s\in\R$ with $|t-s|<\delta$, we have
\Eq{ts}{
\sum_{k=0}^n d_\Z(2^kt)S_k \subseteq \sum_{k=0}^n d_\Z(2^ks)S_k + [5]V + K.
}

By the assumption (ii), there exists a positive integer $m$ such that, for $n\geq m$, \eq{iii} holds.
The members of the sequence $(S_k)$ are closedly $K$-lower bounded, and the family of closedly $K$-lower
bounded sets is closed under finite union, under algebraic sum and under multiplication by scalars, hence there exists 
a bounded set $H\subseteq Y$ such that 
\Eq{*}{
  \bigcup_{\ell=0}^{m-1}\sum_{k=0}^{\ell} 2^kS_k \subseteq \cl(H+K)\subseteq H+V+K.
} 
For the bounded set $H$ there exists a positive number $\delta\leq 1$ such that $\delta H\subseteq V$.
Thus,
\Eq{dd}{
  \delta\bigg(\bigcup_{\ell=0}^{m-1}\sum_{k=0}^{\ell} 2^kS_k\bigg) \subseteq \delta(H+V+K) \subseteq [2]V+K.
}
In view of the Lipschitz property of the distance function $d_\Z$ and by the inequality $0\leq d_Z\leq \frac12$, 
for all $p,q\in\R$, we have that
\Eq{pq}{
  d_\Z(p) \leq d_\Z(q)+|d_\Z(p)-d_\Z(q)| \leq d_\Z(q) +\min\big\{1,|p-q|\big\}.
}
For a closedly $K$-starshaped set $S$, for $0\leq c\leq d$ and for any $W\in\U(Y)$, we have
\Eq{*}{
  cS\subseteq dS+W+K.
}
Indeed, this inclusion is obvious for $c=d$. If $c<d$, then, 
\Eq{*}{
  cS=d\Big(\frac{c}{d}S\Big)\subseteq d\cl(S+K) \subseteq d\Big(S+\frac1d W+K\Big)=dS+W+K.
}
Let $n\in\N$ and $0\leq k\leq n$. Choose $W\in\U(Y)$ such that $[n+1]W\subseteq V$. 
Using inequality \eq{pq} and that the set $S_k$ is closedly $K$-starshaped, for all $k\geq0$ ad for all 
$s,t\in\R$, we obtain
\Eq{*}{
 d_\Z(2^kt)S_k &\subseteq \big(d_\Z(2^ks) +\min\big\{1,2^k|t-s|\}\big)S_k + W+K \\
               &\subseteq d_\Z(2^ks)S_k +\min\big\{1,2^k|t-s|\big\}S_k +W+K.
}
Hence, 
\Eq{split}{
  \sum_{k=0}^n d_\Z(2^kt)S_k  
     \subseteq \sum_{k=0}^n d_\Z(2^ks)S_k  + \sum_{k=0}^n\min\big\{1,2^k|t-s|\big\}S_k + V+K.
}
To complete the proof of \eq{ts}, it suffices to show that, for all $n\geq0$,
\Eq{fs}{
  \sum_{k=0}^n\min\big\{1,2^k|t-s|\big\}S_k\subseteq [4]V+K
}
whenever $|t-s|<\delta$. 

Let $t,s\in\R$ be such that $|t-s|\leq\delta$. We have that $\min\big\{1,2^k|t-s|\big\}\leq 1$ and
$\min\big\{1,2^k|t-s|\big\}\leq 2^{k}\delta$ for all $k\geq0$.  
Using that $S_k$ is closedly $K$-starshaped again and the estimates \eq{iii} and \eq{dd}, we get
\Eq{*}{
  \sum_{k=0}^n\min\big\{1,2^k|t-s|\big\}S_k  
    &\subseteq\sum_{k=0}^{\min\{n,m-1\}}\big(2^{k}\delta S_k+W+K\big)  
       + \sum_{k=\min(n+1,m)}^n\big(S_k+W+K\big) \\
    &\subseteq \delta\sum_{k=0}^{\min\{n,m-1\}} 2^kS_k +[2]V+K \\
    &\subseteq \delta\bigg(\bigcup_{\ell=0}^{m-1}\sum_{k=0}^{\ell} 2^kS_k\bigg) +[2]V+K\subseteq [4]V+K.
} 
This completes the proof of \eq{fs} and hence \eq{ts} holds for all $n\geq0$. Applying \eq{ts}, for all $n\geq0$, it 
follows that
\Eq{*}{
  \sum_{k=0}^n d_\Z(2^kt)S_k 
  \subseteq \cl\bigg(\bigcup_{\ell=0}^\infty\sum_{k=0}^\ell d_\Z(2^ks)S_k\bigg) + [5]V + K.
}
Therefore,
\Eq{*}{
  \cl\bigg(\bigcup_{n=0}^\infty\sum_{k=0}^n d_\Z(2^kt)S_k\bigg)
  \subseteq \bigcup_{n=0}^\infty\sum_{k=0}^n d_\Z(2^kt)S_k + V
  \subseteq \cl\bigg(\bigcup_{\ell=0}^\infty\sum_{k=0}^\ell d_\Z(2^ks)S_k\bigg) + U + K.
}
Thus the proof of inclusion \eq{ts0} is complete.
\end{proof}

In the following lemma we describe an important case when the serial $K$-Cauchy property of a sequence of sets can 
be established.

\Lem{KC}{Assume that $Y$ is a locally convex topological vector space. Let $K\subseteq Y$ be a convex cone and 
$S\subseteq Y$ be a closedly $K$-lower bounded set. Let $(\lambda_n)$ be a sequence of 
nonnegative numbers such that the series $\sum \lambda_n$ is convergent. Define $S_n:=\lambda_n S+K$ for $n\geq 0$. Then 
the sequence $(S_n)$ is serially $K$-Cauchy.}

\begin{proof} By the boundedness assumption on $S$, there exists a bounded set $H\subseteq Y$ such that 
$S\subseteq\cl(H+K)$.

To prove the $K$-Cauchy property of $(S_n)$, fix a neighborhood $U\in\U(Y)$ arbitrarily. Then choose a convex 
$V\in\U(Y)$ such that $[2]V\subseteq U$ and $t>0$ such that $tH\subseteq V$. By the convergence of the series $\sum 
\lambda_n$, there
exists $m\in\N$ such that, for all $n\geq m$,
\Eq{*}{
  t_n:=\sum_{k=m}^n \lambda_k<t.
}
Then, for all $n\geq m$,
\Eq{*}{
 \sum_{k=m}^n S_k
   &=\sum_{k=m}^n (\lambda_k S) + K
   \subseteq \sum_{k=m}^n (\lambda_k \cl(H+K)) + K
   \subseteq \cl\bigg(\sum_{k=m}^n (\lambda_k H) + K\bigg)\\
   &\subseteq \sum_{k=m}^n (\lambda_k H) + V + K
   \subseteq \sum_{k=m}^n \Big(\frac{\lambda_k}{t}V\Big) + \frac{t-t_n}{t}\{0\} + V + K
   \subseteq [2]V+K\subseteq U+K.
}
This shows that \eq{iii} holds and hence $(S_k)$ is serially $K$-Cauchy.
\end{proof}

\section{Diadic convexity and concavity}
\setcounter{theorem}{0}

The first pair of main results of this paper are contained in the following two theorems. The notation $\D$ stands for 
the set of diadic rational numbers, i.e., numbers of the form $\frac{k}{2^n}$, where $k\in\Z$ and $n\in\N$.

\Thm{ConvexTab}{Let $D\subseteq X$ be a nonempty convex set and $A,B:(D-D)\to\P_0(Y)$ 
such that the values of the set-valued map $B$ are closedly $K$-convex, where $K:=\overline{\rec}(B)$. 
Let $F:D\to\P_0(Y)$ be a set-valued mapping which satisfies the Jensen-convexity-type inclusion
\Eq{JCV}{
\dfrac{F(x) + F(y)}{2} + A(x-y) \subseteq \cl\Big(F\Big(\dfrac{x+y}{2}\Big) + B(x-y)\Big) 
   \qquad (x,y\in D).
}
Then $F$ satisfies the convexity type inclusion
\Eq{CV}{
 tF(x)&+(1-t)F(y)+\sum_{k=0}^{\infty} 2d_{\Z}\big(2^k t\big)A\Big(\dfrac{x-y}{2^k}\Big) \\
&\subseteq \cl\bigg(F(tx+(1-t)y)+\sum_{k=0}^{\infty} 2d_{\Z}\big(2^k t\big)B\Big(\dfrac{x-y}{2^k}\Big)\bigg)
   \qquad (x,y\in D,\,t\in\D\cap[0,1]).
}
If, in addition, $0\in\cl(A(u)+K)$ for every $u\in D-D$, then
\Eq{CV+}{
 tF(x)+(1-t)F(y)+A^\perp(t,x-y) \subseteq \cl\big(F(tx+(1-t)y)&+B^\perp(t,x-y)\big)\\
   & (x,y\in D,\,t\in\D\cap[0,1]).
}}

\begin{proof} For the proof of \eq{CV}, we are going to show that, for all $x,y\in D$,
and for all integers $n,m$ with $n\geq 1$, $0\leq m\leq 2^n$,
\Eq{CVn}{
 \frac{m}{2^n}F(x)+\Big(1-\frac{m}{2^n}\Big)F(y) 
  &+ \sum_{k=0}^{n-1}2 d_{\Z}\Big(2^k \frac{m}{2^n}\Big)A\Big(\dfrac{x-y}{2^k}\Big) \\
 &\subseteq \cl\bigg(F\Big(\frac{m}{2^n}x+\Big(1-\frac{m}{2^n}\Big)y\Big) 
   + \sum_{k=0}^{n-1}d_{\Z}\Big(2^k \frac{m}{2^n}\Big)B\Big(\dfrac{x-y}{2^k}\Big)\bigg). 
}

Fix $x,y\in D$ arbitrarily. To verify that \eq{CVn} holds, we will proceed by induction on $n$. 
For $m=0$ or for $m=2^n$, the inclusion \eq{CVn} is obvious (because $d_\Z$ vanishes at integer
values). Thus, for $n=1$, we need to check \eq{CVn} only for $m=1$. In that case, \eq{CVn} is
equivalent to the Jensen convexity assumption \eq{JCV}.

Now, suppose that the inclusion \eq{CVn} holds for some $n$ and let us prove that it is also valid for
$n+1$. Let $0<m<2^{n+1}$ be arbitrary. If $m$ is even, i.e., $m=2\ell$ for some $0<\ell<2^n$, then
the statement follows from the inductive assumption, because $\frac{m}{2^{n+1}}=\frac{\ell}{2^n}$.
Therefore, we may assume that $m$ is odd, i.e., $m=2\ell+1$, where $0\leq \ell<2^n$.
Observe that the fraction $\frac{2\ell+1}{2^{n+1}}$ can be represented as
\Eq{am}{
   \frac{2\ell+1}{2^{n+1}}=\frac12\Big(\frac{\ell+1}{2^{n}}+\frac{\ell}{2^{n}}\Big)
}
We prove that, for $k\in\{0,\dots,n-1\}$,
\Eq{JCd}{
   d_\Z\Big(2^k\frac{2\ell+1}{2^{n+1}}\Big)
    =\frac12\Big(d_\Z\Big(2^k\frac{\ell+1}{2^{n}}\Big)+d_\Z\Big(2^k\frac{\ell}{2^{n}}\Big)\Big).
}
In view of \eq{am}, it suffices to show that $d_\Z$ is affine on the interval 
$\Big[2^k\frac{\ell}{2^{n}},2^k\frac{\ell+1}{2^{n}}\Big]$. For this, it is enough to show
that the interior of this interval does not contain an element from $\frac12\Z$. On the contrary,
assume that, for some $j\in\Z$,
\Eq{*}{
   2^k\frac{\ell}{2^{n}}<\frac{j}2<2^k\frac{\ell+1}{2^{n}}.
}
Then we have that $\ell<2^{n-k-1}j<\ell+1$, which is a contradiction because $n-1\geq k$. 

Now, we can start proving that \eq{CVn} holds for $n+1$. Using \eq{am} and \eq{JCd}, we obtain
\Eq{*}{
  \frac{m}{2^{n+1}}F(x)&+\Big(1-\frac{m}{2^{n+1}}\Big)F(y)
   + \sum_{k=0}^{n}2 d_{\Z}\Big(2^k \frac{m}{2^{n+1}}\Big)A\Big(\frac{x-y}{2^k}\Big) \\
  &= \frac{2\ell+1}{2^{n+1}}F(x)+\Big(1-\frac{2\ell+1}{2^{n+1}}\Big)F(y)
   + \sum_{k=0}^{n}2 d_{\Z}\Big(2^k \frac{2\ell+1}{2^{n+1}}\Big)A\Big(\frac{x-y}{2^k}\Big) \\
  &\subseteq \frac12\Big(\frac{\ell+1}{2^{n}}F(x)+\Big(1-\frac{\ell+1}{2^{n}}\Big)F(y)\Big) 
     +\frac12\Big(\frac{\ell}{2^{n}}F(x)+\Big(1-\frac{\ell}{2^{n}}\Big)F(y)\Big)\\
  &\qquad+\sum_{k=0}^{n-1}\Big(d_\Z\Big(2^k\frac{\ell+1}{2^{n}}\Big)+d_\Z\Big(2^k\frac{\ell}{2^{n}}\Big)\Big)
             A\Big(\frac{x-y}{2^k}\Big) + A\Big(\frac{x-y}{2^n}\Big) \\
  &= \frac12\Big(\frac{\ell+1}{2^{n}}F(x)+\Big(1-\frac{\ell+1}{2^{n}}\Big)F(y)
     +\sum_{k=0}^{n-1}2 d_\Z\Big(2^k\frac{\ell+1}{2^{n}}\Big)A\Big(\frac{x-y}{2^k}\Big)\Big)\\
  &\qquad +\frac12\Big(\frac{\ell}{2^{n}}F(x)+\Big(1-\frac{\ell}{2^{n}}\Big)F(y)
   + \sum_{k=0}^{n-1}2d_\Z\Big(2^k\frac{\ell}{2^{n}}\Big)A\Big(\frac{x-y}{2^k}\Big)\Big)
   + A\Big(\frac{x-y}{2^n}\Big). 
}
By applying the inductive assumption \eq{CVn} with $m=\ell+1$ and $m=\ell$, we obtain that
\Eq{JCVa}{
  \frac{m}{2^{n+1}}F(x)&+\Big(1-\frac{m}{2^{n+1}}\Big)F(y)
   + \sum_{k=0}^{n}2 d_{\Z}\Big(2^k \frac{m}{2^{n+1}}\Big)A\Big(\frac{x-y}{2^k}\Big) \\
 &\subseteq \frac12\cl\bigg(F\Big(\frac{\ell+1}{2^{n}}x+\Big(1-\frac{\ell+1}{2^{n}}\Big)y\Big)
     +\sum_{k=0}^{n-1} 2d_\Z\Big(2^k\frac{\ell+1}{2^{n}}\Big)B\Big(\frac{x-y}{2^k}\Big)\bigg)\\
   &\qquad+\frac12\cl\bigg(F\Big(\frac{\ell}{2^{n}}x+\Big(1-\frac{\ell}{2^{n}}\Big)y\Big)
   + \sum_{k=0}^{n-1} 2d_\Z\Big(2^k\frac{\ell}{2^{n}}\Big)B\Big(\frac{x-y}{2^k}\Big)\bigg) 
      + A\Big(\frac{x-y}{2^n}\Big)\\
 &\subseteq \cl\bigg(\frac12\bigg(F\Big(\frac{\ell+1}{2^{n}}x+\Big(1-\frac{\ell+1}{2^{n}}\Big)y\Big)
     +F\Big(\frac{\ell}{2^{n}}x+\Big(1-\frac{\ell}{2^{n}}\Big)y\Big)\bigg)
     +A\Big(\frac{x-y}{2^n}\Big) \\
   &\qquad+ \sum_{k=0}^{n-1}\bigg(d_\Z\Big(2^k\frac{\ell+1}{2^{n}}\Big)B\Big(\frac{x-y}{2^k}\Big)
      +d_\Z\Big(2^k\frac{\ell}{2^{n}}\Big) B\Big(\frac{x-y}{2^k}\Big)\bigg)\bigg). \\
}
Now, using the Jensen convexity assumption \eq{JCV}, it follows that
\Eq{*}{
  \frac12\bigg(F\Big(\frac{\ell+1}{2^{n}}x+\Big(1-\frac{\ell+1}{2^{n}}\Big)y\Big)
     &+F\Big(\frac{\ell}{2^{n}}x+\Big(1-\frac{\ell}{2^{n}}\Big)y\Big)\bigg)
     +A\Big(\frac{x-y}{2^n}\Big) \\
 &\subseteq \cl\bigg(F\Big(\frac{2\ell+1}{2^{n+1}}x+\Big(1-\frac{2\ell+1}{2^{n+1}}\Big)y\Big)
     +B\Big(\frac{x-y}{2^n}\Big)\bigg).
}
On the other hand, using that the values of the set-valued map $B$ are closedly $K$-convex and 
equation \eq{JCd}, we obtain
\Eq{*}{
 \sum_{k=0}^{n-1}&\bigg(d_\Z\Big(2^k\frac{\ell+1}{2^{n}}\Big)B\Big(\frac{x-y}{2^k}\Big)
      +d_\Z\Big(2^k\frac{\ell}{2^{n}}\Big) B\Big(\frac{x-y}{2^k}\Big)\bigg) \\
  &\!\!\!\subseteq \cl\bigg(K+\sum_{k=0}^{n-1}\bigg(d_\Z\Big(2^k\frac{\ell+1}{2^{n}}\Big)
      +d_\Z\Big(2^k\frac{\ell}{2^{n}}\Big)\bigg) B\Big(\frac{x-y}{2^k}\Big)\bigg)
  = \cl\bigg(K+\sum_{k=0}^{n-1}2d_\Z\Big(2^k\frac{2\ell+1}{2^{n+1}}\Big)B\Big(\frac{x-y}{2^k}\Big)\bigg).
}
Combining the above two inclusions with \eq{JCVa} and replacing $2\ell+1$ by $m$, we arrive at
\Eq{*}{
  \frac{m}{2^{n+1}}F(x)&+\Big(1-\frac{m}{2^{n+1}}\Big)F(y)
   + \sum_{k=0}^{n}2 d_{\Z}\Big(2^k \frac{m}{2^{n+1}}\Big)A\Big(\frac{x-y}{2^k}\Big) \\
  &\subseteq   \cl\bigg(F\Big(\frac{m}{2^{n+1}}x+\Big(1-\frac{m}{2^{n+1}}\Big)y\Big)
     +B\Big(\frac{x-y}{2^n}\Big) + K
  +\sum_{k=0}^{n-1}2d_\Z\Big(2^k\frac{m}{2^{n+1}}\Big)B\Big(\frac{x-y}{2^k}\Big)\bigg)\\
  &= \cl\bigg(F\Big(\frac{m}{2^{n+1}}x+\Big(1-\frac{m}{2^{n+1}}\Big)y\Big)
  +\sum_{k=0}^{n}2d_\Z\Big(2^k\frac{m}{2^{n+1}}\Big)B\Big(\frac{x-y}{2^k}\Big)\bigg),
}
which shows that the statement is also valid for $n+1$. This completes the proof of the induction and the first 
assertion of the theorem.

Assume now that, for all $u\in D-D$, we have $0\in\cl(A(u)+K)$. To prove \eq{CV+}, let $x,y\in D$ and 
$t\in[0,1]\cap\D$ be fixed. If $t\in\{0,1\}$, then \eq{CV+} is trivial. in the rest of the proof assume that 
$t\in]0,1[\cap\D$. Then
\Eq{*}{
 tF(x)+(1-t)F(y)&+A^\perp(t,x-y) \\
 &\subseteq \cl\bigg(tF(x)+(1-t)F(y) + \bigcup_{n=0}^{\infty} \sum_{k=0}^{n}
                 2d_{\Z}(2^kt)A\Big(\frac{x-y}{2^k}\Big)\bigg)\\
 &\subseteq \cl\bigg(tF(x)+(1-t)F(y) + \bigcup_{n=0}^{\infty} \sum_{k=0}^{n}
                 2d_{\Z}(2^kt)\cl\Big(A\Big(\frac{x-y}{2^k}\Big)+K\Big)\bigg).
}
Now, since $t$ is a diadic number, there exists an $n\in\N$ and an odd number $\ell$ such that $t=\frac{\ell}{2^{m}}$. 
Then, for $k\geq m$, $d_{\Z}(2^kt)=0$. In addition, $0\in\cl\big(A\big(\frac{x-y)}{2^k}\big)+K\big)$ for 
all $k\in\{0,\dots,m-1\}$, therefore
\Eq{*}{
  \bigcup_{n=0}^{\infty} \sum_{k=0}^{n}2d_{\Z}(2^kt)\cl\Big(A\Big(\frac{x-y}{2^k}\Big)+K\Big)
  &= \sum_{k=0}^{m-1}2d_{\Z}(2^kt)\cl\Big(A\Big(\frac{x-y}{2^k}\Big)+K\Big) \\
  &= \sum_{k=0}^{\infty}2d_{\Z}(2^kt)\cl\Big(A\Big(\frac{x-y}{2^k}\Big)+K\Big).
}
Using this formula, the previous inclusion, the first part of the theorem and $d_\Z(t)>0$ we arrive at
\Eq{AB}{
 tF(x)+(1-t)F(y)&+A^\perp(t,x-y) \\
 &\subseteq \cl\bigg(tF(x)+(1-t)F(y) 
    + \sum_{k=0}^{\infty}2d_{\Z}(2^kt)\cl\Big(A\Big(\frac{x-y}{2^k}\Big)+K\Big)\bigg)\\
 &=\cl\bigg(tF(x)+(1-t)F(y) 
    + \sum_{k=0}^{\infty}2d_{\Z}(2^kt)\Big(A\Big(\frac{x-y}{2^k}\Big)+K\Big)\bigg) \\
 &=\cl\bigg(tF(x)+(1-t)F(y) 
    + \sum_{k=0}^{\infty}2d_{\Z}(2^kt)A\Big(\frac{x-y}{2^k}\Big)+\rec(B)\bigg) \\
 &\subseteq \cl\bigg(F(tx+(1-t)y) 
    + \sum_{k=0}^{\infty}2d_{\Z}(2^kt)B\Big(\frac{x-y}{2^k}\Big)+\rec(B)\bigg) \\
 &= \cl\bigg(F(tx+(1-t)y) 
    + \sum_{k=0}^{\infty}2d_{\Z}(2^kt)B\Big(\frac{x-y}{2^k}\Big)\bigg).
}
On the other hand, we have 
\Eq{*}{
  \sum_{k=0}^{\infty}2d_{\Z}(2^kt)B\Big(\frac{x-y}{2^k}\Big)
  =\sum_{k=0}^{m-1}2d_{\Z}(2^kt)B\Big(\frac{x-y}{2^k}\Big)
  \subseteq \cl\bigg(\bigcup_{n=0}^{\infty} \sum_{k=0}^{n}2d_{\Z}(2^kt)B\Big(\frac{x-y}{2^k}\Big)\bigg)
  =B^\perp(t,x-y),
}
which combined with \eq{AB} implies \eq{CV+} and completes the proof of the theorem.
\end{proof}

To formulate the most useful corollaries of our main results we postulate the following general hypotheses.
\begin{enumerate}[(H1)]
\item $D\subseteq X$ is a nonempty convex set, $K\subseteq Y$ is a nonempty convex cone; 
\item $S_0\subseteq Y$ is a closedly $K$-convex and closedly $K$-starshaped set;
\item $\varphi:(D-D)\to\R_+$ is a nonnegative function such that \eq{phi} holds for all $x\in D-D$. 
\end{enumerate}
Note that, by the convexity of $D$, the set $(D-D)$ is starshaped, thus \prp{Tab} can be applied. 

The next two corollaries are about approximately and strongly $K$-Jensen convex set-valued maps, respectively.

\Cor{Convex+1}{Assume that (H1), (H2), and (H3) hold and $F:D\to\P_0(Y)$ is a set-valued mapping which satisfies
\Eq{*}{
\dfrac{F(x) + F(y)}{2} \subseteq \cl\Big(F\Big(\dfrac{x+y}{2}\Big) 
    + \varphi(x-y)S_0 + K \Big) \qquad (x,y\in D).
}
Then 
\Eq{*}{
 tF(x)+(1-t)F(y) \subseteq \cl\big(F(tx+(1-t)y)+\varphi^\perp(t,x-y)S_0+K\big) \qquad
 (x,y\in D,\,t\in\D\cap[0,1]).
}}

\Cor{Convex+2}{Assume that (H1), (H2), and (H3) hold and $F:D\to\P_0(Y)$ is a set-valued mapping which satisfies
\Eq{*}{
\dfrac{F(x) + F(y)}{2} + \varphi(x-y)S_0 \subseteq \cl\Big(F\Big(\dfrac{x+y}{2}\Big) 
  + K \Big) \qquad (x,y\in D).
}
Then 
\Eq{*}{
 tF(x)+(1-t)F(y) + \varphi^\perp(t,x-y)S_0 \subseteq \cl\big(F(tx+(1-t)y) + K \big)\qquad 
 (x,y\in D,\,t\in\D\cap[0,1]).
}}

\begin{proof}[Proof of the Corollaries \ref{CConvex+1} and \ref{CConvex+2}]
Using the second assertion of \thm{ConvexTab} with the set-valued maps $A(u)=\{0\}$ and $B(u)=\varphi(u)S_0+K$ (resp., 
$A(u)=\varphi(u)S_0$ and $B(u)=K$) and applying the \prp{Tab}, we obtain \cor{Convex+1} 
(resp., \cor{Convex+2}). Note that, in both cases, $0\in\cl(A(u)+K)$ for all $u\in D-D$.
\end{proof}

The next result concerns the case of concavity type inclusions.

\Thm{ConcaveTab}{Let $D\subseteq X$ be a nonempty convex set and $A,B:(D-D)\to\P_0(Y)$ be
such that the values of the set-valued map $B$ are closedly $K$-convex, where $K:=\overline{\rec}(B)$. Let 
$F:D\to\P_0(Y)$ be a set-valued mapping which satisfies the Jensen-concavity-type inclusion
\Eq{JCC}{
F\Big(\dfrac{x+y}{2}\Big) + A(x-y) \subseteq \cl\bigg(\frac{F(x) + F(y)}{2} + B(x-y)\bigg) 
   \qquad (x,y\in D),
}
and has closedly $K$-convex values, i.e, $F(x)$ is closedly $K$-convex, for all $x\in D$.
Then $F$ satisfies the convexity type inclusion
\Eq{CC}{
 F(tx&+(1-t)y)+\sum_{k=0}^{\infty} 2d_{\Z}\big(2^k t\big)A\Big(\dfrac{x-y}{2^k}\Big) \\
&\subseteq \cl\bigg(tF(x)+(1-t)F(y)+\sum_{k=0}^{\infty} 2d_{\Z}\big(2^k t\big)B\Big(\dfrac{x-y}{2^k}\Big)\bigg)
   \qquad (x,y\in D,\,t\in\D\cap[0,1]).
}
If, in addition, $0\in\cl(A(u)+K)$ for every $u\in D-D$, then
\Eq{CC+}{
 F(tx+(1-t)y)+A^\perp(t,x-y) \subseteq \cl\big(tF(x)+(1-t)F(y)&+B^\perp(t,x-y)\big)\\
   & (x,y\in D,\,t\in\D\cap[0,1]).
}}

\begin{proof}
For the proof of \eq{CC}, we are going to show that, for all $x,y\in D$,
and for all integers $n,m$ with $n\geq 1$, $0\leq m\leq 2^n$,
\Eq{CCn}{
  F\Big(\frac{m}{2^n}x+\Big(1-\frac{m}{2^n}\Big)y\Big)
  &+ \sum_{k=0}^{n-1}2 d_{\Z}\Big(2^k \frac{m}{2^n}\Big)A\Big(\dfrac{x-y}{2^k}\Big) \\
 &\subseteq \cl\bigg(\frac{m}{2^n}F(x)+\Big(1-\frac{m}{2^n}\Big)F(y) 
   + \sum_{k=0}^{n-1}2d_{\Z}\Big(2^k \frac{m}{2^n}\Big)B\Big(\dfrac{x-y}{2^k}\Big)\bigg). 
}
Fix $x,y\in D$ arbitrarily. To verify that \eq{CCn} holds, we will proceed by induction on $n$. For 
$n=1$ we have that $0\leq m\leq 2$, but if $m=0$ or $m=2$ then, equation \eq{CCn} follows immediately. 
Thus we have only to check that \eq{CCn} is valid for $m=1$, which is trivial because for $n=1$ and $m=1$
\eq{CCn} is the same as \eq{JCC}. Now suppose that \eq{CCn} holds for $n\geq1$ and $0\leq m \leq 2^n$, and 
let us prove that it is also valid for $n+1$ and for $0\leq m\leq 2^{n+1}$. As in the proof for the 
convexity type inclusion, it will be enough to consider the case when $m$ has the form $m=2\ell +1$ for some 
$\ell\in\N\cup\{0\}$. Now we can start our proof for $n+1$ using the relations \eq{am} and \eq{JCd}, to obtain 
\Eq{*}{
 F\Big(\frac{2\ell +1}{2^{n+1}}x&+\Big(1-\frac{2\ell +1}{2^{n+1}}\Big)y\Big)
  + \sum_{k=0}^{n}2 d_{\Z}\Big(2^k \frac{2\ell +1}{2^{n+1}}\Big)A\Big(\dfrac{x-y}{2^k}\Big) \\ 
&=F\Bigg(\frac12\bigg(\Big(\frac{\ell}{2^{n}} + \frac{\ell +1}{2^{n}}\Big)x
	+	\Big(2-\frac{\ell +1}{2^{n}}-\frac{\ell}{2^{n}}\Big)y\bigg)\Bigg)
  + \sum_{k=0}^{n}2 d_{\Z}\Big(2^k \frac{2\ell +1}{2^{n+1}}\Big)A\Big(\dfrac{x-y}{2^k}\Big) \\
&=F\Bigg(\frac12\bigg[\frac{\ell}{2^{n}}x + \Big(1-\frac{\ell}{2^{n}}\Big)y
	+	\frac{\ell +1}{2^{n}}x+\Big(1-\frac{\ell +1}{2^{n}}\Big)y\bigg]\Bigg) + A\Big(\frac{x-y}{2^n}\Big) \\
&\qquad+ \sum_{k=0}^{n-1}\bigg[d_\Z\Big(2^k\frac{\ell+1}{2^{n}}\Big)+d_\Z\Big(2^k\frac{\ell}{2^{n}}\Big)\bigg]
	A\Big(\dfrac{x-y}{2^k}\Big) \\
&\subseteq F\Bigg(\frac12\bigg[\frac{\ell}{2^{n}}x + \Big(1-\frac{\ell}{2^{n}}\Big)y
	+	\frac{\ell +1}{2^{n}}x+\Big(1-\frac{\ell +1}{2^{n}}\Big)y\bigg]\Bigg) + A\Big(\frac{x-y}{2^n}\Big) \\
&\qquad+ \sum_{k=0}^{n-1}d_\Z\Big(2^k\frac{\ell}{2^{n}}\Big)A\Big(\dfrac{x-y}{2^k}\Big)
	+ \sum_{k=0}^{n-1}d_\Z\Big(2^k\frac{\ell+1}{2^{n}}\Big)A\Big(\dfrac{x-y}{2^k}\Big).
}
By the Jensen concavity property of $F$, we have that 
\Eq{JCCp}{
F\Bigg(\frac12\bigg[\frac{\ell}{2^{n}}x &+ \Big(1-\frac{\ell}{2^{n}}\Big)y
	+\frac{\ell +1}{2^{n}}x+\Big(1-\frac{\ell +1}{2^{n}}\Big)y\bigg]\Bigg) + A\Big(\frac{x-y}{2^n}\Big) \\
&\subseteq \cl\Bigg(\frac12 F\Big( \frac{\ell}{2^{n}}x + \Big(1-\frac{\ell}{2^{n}}\Big)y \Big)
	+\frac12 F\Big( \frac{\ell+1}{2^{n}}x+\Big(1-\frac{\ell +1}{2^{n}}\Big)y \Big)
	+ B\Big(\frac{x-y}{2^n}\Big)\Bigg).
}
Therefore, using \eq{JCCp}, we can obtain the following inclusions
\Eq{*}{
F\Big(\frac{2\ell +1}{2^{n+1}}x&+\Big(1-\frac{2\ell +1}{2^{n+1}}\Big)y\Big)
  + \sum_{k=0}^{n}2 d_{\Z}\Big(2^k \frac{2\ell +1}{2^{n+1}}\Big)A\Big(\dfrac{x-y}{2^k}\Big) \\ 
&\subseteq \cl\Bigg(\frac12 F\Big( \frac{\ell}{2^{n}}x + \Big(1-\frac{\ell}{2^{n}}\Big)y \Big)
  +\frac12 F\Big( \frac{\ell +1}{2^{n}}x+\Big(1-\frac{\ell +1}{2^{n}}\Big)y \Big)+ B\Big(\frac{x-y}{2^n}\Big)\Bigg) \\
&\qquad+ \sum_{k=0}^{n-1}d_\Z\Big(2^k\frac{\ell}{2^{n}}\Big)A\Big(\dfrac{x-y}{2^k}\Big)
	+ \sum_{k=0}^{n-1}d_\Z\Big(2^k\frac{\ell+1}{2^{n}}\Big)A\Big(\dfrac{x-y}{2^k}\Big)	\\
&\subseteq \cl\Bigg( \frac12\bigg[ F\Big(\frac{\ell}{2^{n}}x + \Big(1-\frac{\ell}{2^{n}}\Big)y \Big)
	+\sum_{k=0}^{n-1}2d_\Z\Big(2^k\frac{\ell}{2^{n}}\Big)A\Big(\dfrac{x-y}{2^k}\Big) \\
&\qquad + F\Big( \frac{\ell +1}{2^{n}}x+\Big(1-\frac{\ell +1}{2^{n}}\Big)y \Big)
	+\sum_{k=0}^{n-1}2d_\Z\Big(2^k\frac{\ell+1}{2^{n}}\Big)A\Big(\dfrac{x-y}{2^k}\Big)\bigg]
	+B\Big(\frac{x-y}{2^n}\Big)\Bigg).
}
Our inductive assumption for $m=\ell$ and $m=\ell+1$ gives us the following relations,
\Eq{Il}{
F\Big( \frac{\ell}{2^{n}}x + \Big(1-\frac{\ell}{2^{n}}\Big)y \Big)
	&+ \sum_{k=0}^{n-1}2d_\Z\Big(2^k\frac{\ell}{2^{n}}\Big)A\Big(\dfrac{x-y}{2^k}\Big)\\
&\subseteq \cl\Bigg(\frac{\ell}{2^{n}}F(x) +  \Big(1-\frac{\ell}{2^{n}}\Big)F(y)
	+ \sum_{k=0}^{n-1}2d_\Z\Big(2^k\frac{\ell}{2^{n}}\Big)B\Big(\dfrac{x-y}{2^k}\Big)\Bigg),
}
and
\Eq{Il+}{
F\Big( \frac{\ell+1}{2^{n}}x &+ \Big(1-\frac{\ell+1}{2^{n}}\Big)y \Big)
	+ \sum_{k=0}^{n-1}2d_\Z\Big(2^k\frac{\ell+1}{2^{n}}\Big)A\Big(\dfrac{x-y}{2^k}\Big)\\
&\subseteq \cl\Bigg(\frac{\ell+1}{2^{n}}F(x) +  \Big(1-\frac{\ell+1}{2^{n}}\Big)F(y)
	+ \sum_{k=0}^{n-1}2d_\Z\Big(2^k\frac{\ell+1}{2^{n}}\Big)B\Big(\dfrac{x-y}{2^k}\Big)\Bigg).
}
Thus, by \eq{JCd}, by using that $F$ and $B$ have closedly $K$-convex values, and by  
\eq{Il} and \eq{Il+}, we arrive at
\Eq{*}{
 F&\Big(\frac{2\ell +1}{2^{n+1}}x+\Big(1-\frac{2\ell +1}{2^{n+1}}\Big)y\Big)
  + \sum_{k=0}^{n}2 d_{\Z}\Big(2^k \frac{2\ell +1}{2^{n+1}}\Big)A\Big(\dfrac{x-y}{2^k}\Big) \\
&\subseteq \cl\Bigg( \frac12\Bigg[\frac{\ell}{2^{n}}F(x) +  \Big(1-\frac{\ell}{2^{n}}\Big)F(y)
	+ \sum_{k=0}^{n-1}2d_\Z\Big(2^k\frac{\ell}{2^{n}}\Big)B\Big(\dfrac{x-y}{2^k}\Big) \\
&\qquad + \frac{\ell+1}{2^{n}}F(x) +  \Big(1-\frac{\ell+1}{2^{n}}\Big)F(y)
	+ \sum_{k=0}^{n-1}2d_\Z\Big(2^k\frac{\ell+1}{2^{n}}\Big)B\Big(\dfrac{x-y}{2^k}\Big)\Bigg]
	+ B\Big(\frac{x-y}{2^n}\Big)\Bigg) \\
&\subseteq \cl\Bigg( \frac12\Bigg[\frac{2\ell+1}{2^{n}}F(x) +  \Big(2-\frac{2\ell+1}{2^{n}}\Big)F(y)\\ 
&\qquad + \sum_{k=0}^{n-1}2\Big(d_\Z\Big(2^k\frac{\ell}{2^{n}}\Big) + d_\Z\Big(2^k\frac{\ell+1}{2^{n}}\Big)
	  \Big)B\Big(\dfrac{x-y}{2^k}\Big) \Bigg] + K + B\Big(\frac{x-y}{2^n}\Big) \Bigg) \\
&= \cl\Bigg( \frac{2\ell+1}{2^{n+1}}F(x) +  \Big(1-\frac{2\ell+1}{2^{n+1}}\Big)F(y) 
	+ \sum_{k=0}^{n}2d_\Z\Big(2^k\frac{2\ell+1}{2^{n+1}}\Big) B\Big(\dfrac{x-y}{2^k}\Big)\Bigg),
}
which shows that the statement is also valid for $n+1$. Thus the induction, and therefore, 
the proof of the theorem is complete.

The last statement of the theorem can be proved in the same manner as that of \thm{ConvexTab}.
\end{proof}

The next two corollaries are about approximately and strongly $K$-Jensen concave set-valued mapping, 
respectively.

\Cor{Concave+1}{Assume that (H1), (H2), and (H3) hold and let $F:D\to\P_0(Y)$ be a set-valued mapping with 
closedly $K$-convex values satisfying
\Eq{*}{
F\Big(\dfrac{x+y}{2}\Big) \subseteq \cl\Big(\dfrac{F(x) + F(y)}{2} 
  + \varphi(x-y)S_0 + K \Big) \qquad (x,y\in D).
}
Then 
\Eq{*}{
 F(tx+(1-t)y) \subseteq \cl\big(tF(x)+(1-t)F(y) + \varphi^\perp(t,x-y)S_0 + K \big)\qquad  
 (x,y\in D,\,t\in\D\cap[0,1]).
}}

\Cor{Concave+2}{Assume that (H1), (H2), and (H3) hold and $F:D\to\P_0(Y)$ be a set-valued mapping with 
closedly $K$-convex values satisfying
\Eq{*}{
F\Big(\dfrac{x+y}{2}\Big) + \varphi(x-y)S_0 \subseteq \cl\Big(\dfrac{F(x) + F(y)}{2} 
    + K\Big) \qquad (x,y\in D).
}
Then 
\Eq{*}{
 F(tx+(1-t)y)+\varphi^\perp(t,x-y)S_0 \subseteq \cl\big(tF(x)+(1-t)F(y) + K \big)\qquad  
 (x,y\in D,\,t\in\D\cap[0,1]).
}}

\begin{proof}[Proof of the Corollaries \ref{CConcave+1} and \ref{CConcave+2}]
Using \thm{ConcaveTab} with the set-valued maps $A(u)=\{0\}$ and $B(u)=\varphi(u)S_0+K$ (resp., 
$A(u)=\varphi(u)S_0$ and $B(u)=K$) and applying the \prp{Tab}, we obtain \cor{Concave+1} 
(resp., \cor{Concave+2}). In both settings, we have that $K\subseteq\overline{\rec}(B)$, thus
the assumption that the values of $F$ are closedly $K$-convex implies that the values are also closedly 
$\overline{\rec}(B)$-convex.
\end{proof}

\section{Bernstein--Doetsch type results}
\setcounter{theorem}{0}

The two result below deduce convexity and concavity properties of the set-valued map $F$ based on \thm{ConvexTab} and 
\thm{ConcaveTab}.

\Thm{ConvexTab+}{
Let $D\subseteq X$ be a nonempty convex set, $A,B:(D-D)\to\P_0(Y)$ and denote $\overline{\rec}(B)$
by $K$. Let $F:D\to\P_0(Y)$ be a set-valued mapping which satisfies the Jensen-convexity-type inclusion
\Eq{JCV+}{
\dfrac{F(x) + F(y)}{2} + A(x-y) \subseteq \cl\Big(F\Big(\dfrac{x+y}{2}\Big) + B(x-y)\Big) 
   \qquad (x,y\in D).
}
Assume in addition that the maps $A,B$ and $F$ have the following properties
\begin{enumerate}[(i)]
\item For all $x\in D$, $F(x)$ is closedly $K$-starshaped with respect to some 
 element of $Y$ and also closedly $K$-lower bounded;
\item $F$ is directionally $K$-upper semicontinuous on $D$;
\item For all $u\in D-D$, the sets $A(u)$ and $B(u)$ are closedly $K$-starshaped, closedly-$K$-lower bounded and 
 $B(u)$ is also closedly $K$-convex;
\item For all $u\in D-D$, the sequences $\Big(A\Big(\dfrac{u}{2^k}\Big)\Big)$ and 
 $\Big(B\Big(\dfrac{u}{2^k}\Big)\Big)$ are serially $K$-Cauchy.
\end{enumerate}
Then $F$ satisfies the following convexity type inclusion
\Eq{JCAB}{
tF(x) + (1-t)F(y)+ A^\perp(t,x-y) \subseteq \cl\big(F(tx+(1-t)y)+B^\perp(t,x-y)+K\big) \qquad \\
(x,y\in D,\,t\in[0,1]).
}}

\begin{proof}
Let $x,y\in D$, $t\in[0,1]$ be fixed. Let $U\in\U(Y)$ be arbitrary and choose $V\in\U(Y)$ such that 
$[5]V\subseteq U$. By \lem{dusc}, there exists $\delta_1\in]0,t]$ such that, for $s\in[t-\delta_1,t]$,
\Eq{tsF}{
tF(x) + (1-t)F(y)\subseteq sF(x)+(1-s)F(y) + V + K.
}
Applying \lem{CchyConv} for the sequences $S_k:=A\big(\frac{u}{2^k}\big)$ and $S_k:=B\big(\frac{u}{2^k}\big)$, 
there exists $\delta_2\in]0,\delta_1]$ such that, for $s\in[t-\delta_2,t]$,
\Eq{tsTA}{
  A^\perp(t,x-y)\subseteq A^\perp(s,x-y)+V+K \qquad\mbox{and}\qquad B^\perp(s,x-y)\subseteq B^\perp(t,x-y)+V+K.
}
Adding up \eq{tsF} and the first inclusion in \eq{tsTA} side by side, for $s\in[t-\delta_2,t]$, we obtain
\Eq{*}{
  tF(x) + (1-t)F(y) + A^\perp(t,x-y) \subseteq sF(x)+(1-s)F(y) + A^\perp(s,x-y)+[2]V+K.
}
As, for $s\in\D\cap[0,1]$, by \thm{ConvexTab}, we have that 
\Eq{*}{
sF(x)+(1-s)F(y)+A^\perp(s,x-y) \subseteq \cl\big(F(sx+(1-s)y)+B^\perp(s,x-y)\big).
}
Therefore, for $s\in[t-\delta_2,t]\cap\D$, we get
\Eq{tAsB}{
tF(x) + (1-t)F(y) + A^\perp(t,x-y) \subseteq F(sx+(1-s)y) + B^\perp(s,x-y)+[3]V+K.
}
Observing that $F(sx+(1-s)y)=F(tx+(1-t)y+(s-t)(x-y))$ and then using the  $K$-upper semicontinuity of the set-valued map 
$F$ at the point $tx+(1-t)y$ along the direction $x-y$, there exists $\delta_3\in]0,\delta_2]$ such that, for 
$s\in[t-\delta_3,t]$,
\Eq{FKuhc}{
F(sx+(1-s)y)\subseteq F(tx+(1-t)y) + V+K.
}
Combining this inclusion and the second one in \eq{tsTA} with \eq{tAsB} and chosing $s\in[t-\delta_3,t]\cap\D$ 
arbitrarily, we get
\Eq{*}{
  tF(x) + (1-t)F(y) + A^\perp(t,x-y) 
  &\subseteq F(tx+(1-t)y) + B^\perp(t,x-y)+[5]V+K \\
  &\subseteq F(tx+(1-t)y) + B^\perp(t,x-y)+U.
}
Since $U$ was arbitrary, the above inclusion implies that \eq{JCAB} holds.
\end{proof}

To obtain some corollaries of our Bernstein--Doetsch type theorems of this section, we keep hypothesis (H3) 
that we postulated after \thm{ConvexTab} and replace (H1) and (H2) by the following stronger requirements.
\begin{enumerate}
\item[(H1$^*$)] $D\subseteq X$ is a nonempty convex set, $Y$ is a locally convex topological linear space and 
$K\subseteq Y$ is a nonempty convex cone; 
\item[(H2$^*$)] $S_0\subseteq Y$ is a closedly $K$-convex, closedly $K$-starshaped, and closedly $K$-lower bounded 
set.
\end{enumerate}

\Cor{Convex+1+}{Assume that (H1$^*$), (H2$^*$), and (H3) hold and $F:D\to\P_0(Y)$ is a directionally $K$-upper 
semicontinuous set-valued map such that, for all $x\in D$, that $F(x)$ is closedly $K$-starshaped with respect to some 
element of $Y$ and also closedly $K$-lower bounded. If $F$ satisfies
\Eq{*}{
\dfrac{F(x) + F(y)}{2} \subseteq \cl\Big(F\Big(\dfrac{x+y}{2}\Big) 
    + \varphi(x-y)S_0 + K \Big) \qquad (x,y\in D),
}
then 
\Eq{*}{
 tF(x)+(1-t)F(y) \subseteq \cl\big(F(tx+(1-t)y)+\varphi^\perp(t,x-y)S_0+K\big) \qquad
 (x,y\in D,\,t\in[0,1]).
}}

\Cor{Convex+2+}{Assume that (H1$^*$), (H2$^*$), and (H3) hold and $F:D\to\P_0(Y)$ is a directionally $K$-upper 
semicontinuous set-valued map such that, for all $x\in D$, that $F(x)$ is closedly $K$-starshaped with respect to some 
element of $Y$ and also closedly $K$-lower bounded. If $F$ satisfies
\Eq{*}{
\dfrac{F(x) + F(y)}{2} + \varphi(x-y)S_0 \subseteq \cl\Big(F\Big(\dfrac{x+y}{2}\Big) 
  + K \Big) \qquad (x,y\in D),
}
then 
\Eq{*}{
 tF(x)+(1-t)F(y) + \varphi^\perp(t,x-y)S_0 \subseteq \cl\big(F(tx+(1-t)y) + K \big)\qquad 
 (x,y\in D,\,t\in[0,1]).
}}

\begin{proof}[Proof of the Corollaries \ref{CConvex+1} and \ref{CConvex+2}]
We are going to apply \thm{ConvexTab+} with the set-valued maps $A(u):=\{0\}$ and $B(u):=\varphi(u)S_0+K$ 
(resp., $A(u):=\varphi(u)S_0$ and $B(u):=K$). Then, conditions (i) and (ii) of this theorem are assumed for $F$ and 
(iii) trivially follows from (H2$^*$) and (H3$^*$), respectively.
For the proof of condition (iv), we apply \lem{KC} for the set $S:=\{0\}$, $S:=S_0$, $S:=S_0$, and $S:=\{0\}$, 
respectively, and for the sequence $\lambda_k:=\varphi\big(\frac{x-y}{2^k}\big)$ (in each case) and use that $\varphi$ 
satisfies the convergence condition \eq{phi} by hypothesis (H4). 
 
Thus, by \thm{ConvexTab+}, \eq{JCV+} implies \eq{JCAB}, which, applying the \prp{Tab}, yields the corresponding 
convexity type inclusion in \cor{Convex+1} (resp., \cor{Convex+2}). Note that, in both cases, $0\in\cl(A(u)+K)$ for all 
$u\in D-D$.
\end{proof}

\Thm{ConcaveTab+}{Let $D\subseteq X$ be a nonempty convex set and $A,B:(D-D)\to\P_0(Y)$ and
denote $\overline{\rec}(B)$ by $K$. Let $F:D\to\P_0(Y)$ be a set-valued mapping which
satisfies the Jensen-concavity-type inclusion
\Eq{JCCV}{
F\bigg(\dfrac{x+y}{2}\bigg)+ A(x-y) \subseteq 
\cl\bigg(\dfrac{F(x) + F(y)}{2}  + B(x-y)\bigg) \qquad (x,y\in D).
}
Assume in addition that the maps $A,B$ and $F$ have the following properties
\begin{enumerate}[(i)]
\item For all $x\in D$, $F(x)$ is closedly $K$-convex and closedly $K$-starshaped with respect to some 
 element of $Y$ and also closedly $K$-lower bounded;
\item $F$ is directionally $K$-upper semicontinuous on $D$;
\item For all $u\in D-D$, the sets $A(u)$ and $B(u)$ are closedly $K$-starshaped, 
 closedly-$K$-lower bounded and $B(u)$ is also closedly $K$-convex;
\item For all $u\in D-D$, the sequences $\Big(A\Big(\dfrac{u}{2^k}\Big)\Big)$ and 
 $\Big(B\Big(\dfrac{u}{2^k}\Big)\Big)$ are serially $K$-Cauchy.
\end{enumerate}
Then $F$ satisfies the following concavity type inclusion
\Eq{JCCAB}{
F(tx+(1-t)y)+ A^\perp(t,x-y) \subseteq\cl\big(tF(x) + (1-t)F(y)+B^\perp(t,x-y)+K\big) \qquad \\
(x,y\in D,\,t\in[0,1]).
}}

\begin{proof} Using \thm{ConcaveTab} and \lem{dusc} and \lem{CchyConv}, the proof of the theorem is completely 
analogous to that of \thm{ConvexTab+}. 
\end{proof}

\Cor{Concave+1+}{Assume that (H1$^*$), (H2$^*$), and (H3) hold and $F:D\to\P_0(Y)$ is a directionally $K$-upper 
semicontinuous set-valued map such that, for all $x\in D$, that $F(x)$ is closedly $K$-convex, closedly $K$-starshaped 
with respect to some element of $Y$ and also closedly $K$-lower bounded. If $F$ satisfies
\Eq{*}{
F\Big(\dfrac{x+y}{2}\Big) \subseteq \cl\Big(\dfrac{F(x) + F(y)}{2} 
  + \varphi(x-y)S_0 + K \Big) \qquad (x,y\in D)
}
then 
\Eq{*}{
 F(tx+(1-t)y) \subseteq \cl\big(tF(x)+(1-t)F(y) + \varphi^\perp(t,x-y)S_0 + K \big)\qquad  
 (x,y\in D,\,t\in[0,1]).
}}

\Cor{Concave+2+}{Assume that (H1$^*$), (H2$^*$), and (H3) hold and $F:D\to\P_0(Y)$ is a directionally $K$-upper 
semicontinuous set-valued map such that, for all $x\in D$, that $F(x)$ is closedly $K$-convex, closedly $K$-starshaped 
with respect to some element of $Y$ and also closedly $K$-lower bounded. If $F$ satisfies
\Eq{*}{
F\Big(\dfrac{x+y}{2}\Big) + \varphi(x-y)S_0 \subseteq \cl\Big(\dfrac{F(x) + F(y)}{2} 
    + K\Big) \qquad (x,y\in D),
}
then 
\Eq{*}{
 F(tx+(1-t)y)+\varphi^\perp(t,x-y)S_0 \subseteq \cl\big(tF(x)+(1-t)F(y) + K \big)\qquad  
 (x,y\in D,\,t\in[0,1]).
}}

\begin{proof}[Proof of the Corollaries \ref{CConcave+1} and \ref{CConcave+2}]
We apply \thm{ConcaveTab+} with the set-valued maps $A(u):=\{0\}$ and $B(u):=\varphi(u)S_0+K$ 
(resp., $A(u):=\varphi(u)S_0$ and $B(u):=K$). Argueing as in the proof of Corollaries \ref{CConvex+1+} 
and \ref{CConvex+2+}, we get that the conditions of \thm{ConcaveTab+} are satisfied in each setting.

Thus, by \thm{ConvexTab+}, \eq{JCCV} implies \eq{JCCAB}, which, applying the \prp{Tab}, we obtain 
the expected concavity type inclusions. In both settings, we have that $K\subseteq\overline{\rec}(B)$, thus
the assumption that the values of $F$ are closedly $K$-convex implies that the values are also closedly 
$\overline{\rec}(B)$-convex.
\end{proof}


\def\cprime{$'$} 
\providecommand{\bysame}{\leavevmode\hbox to3em{\hrulefill}\thinspace}
\def\MR#1{}
\providecommand{\MRhref}[2]{%
  \href{http://www.ams.org/mathscinet-getitem?mr=#1}{#2}
}
\providecommand{\href}[2]{#2}

\end{document}